\documentclass[12pt]{article}
\usepackage{amssymb, amsmath, amsfonts,dsfont, mathrsfs,euscript,eufrak,colonequals}
\usepackage[T1]{fontenc}
\usepackage{cmap}
\usepackage{textcomp}
\usepackage{latexsym}
\usepackage{indentfirst}

\newtheorem{Theorem}{Theorem}

\newtheorem{Remark}{Remark}

\newcommand{\N}{\mathbb N}
\newcommand{\R}{\mathbb R}
\newcommand{\Z}{\mathbb Z}

\newcommand{\CC}{\mathbb C}

\newcommand{\defeq}{\colonequals}
\newcommand{\dd}{d}
\newcommand{\Ln}{\mathrm{Ln}\,}
\newcommand{\Arg}{\mathrm{Arg}\,}

\newcommand{\tovar}{\xrightarrow{\mathrm{var}}}

\newcommand{\e}{\varepsilon}
\newcommand{\la}{\langle}
\newcommand{\ra}{\rangle}
\newcommand{\D}{\boldsymbol{D}}
\newcommand{\DS}{\boldsymbol{D}_{\!S}}

\newcommand{\UnionXn}{\mathcal{U}}

\begin{document}
	
	\author{I. A. Alexeev\footnote{The Euler International Mathematical Institute, 10 Pesochnaya nab., 197022 St. Petersburg, Russia, e-mail: \texttt{vanyalexeev@list.ru}},\qquad A. A. Khartov\footnote{Smolensk State University, 4 Przhevalsky st., 214000 Smolensk, Russia, e-mail: \texttt{alexeykhartov@gmail.com} }}

\title{Spectral representations of characteristic functions of discrete probability laws}

\maketitle
		
		\begin{abstract}
			We consider discrete probability laws on the real line, whose characteristic functions are separated from zero. This class includes arbitrary discrete infinitely divisible laws and lattice probability laws, whose characteristic functions have no zeroes on the real line. We show that characteristic functions of such laws admit spectral L\'evy--Khinchine type representation with non-monotonic L\'evy spectral function. We also apply the representations of such laws to obtain limit and compactness theorems  with convergence in variation to probability laws  from this class.
		\end{abstract}
		

\section{Introduction}

This paper is devoted to the L\'evy--Khinchine type representations for characteristic functions of univariate discrete probability laws from rather wide class.

Let $F$ be a  distribution function of a probability law on the real line $\R$. Recall that $F$ is called \textit{infinitely divisible} if for every positive integer $n$ there exists a  distribution function $F_{1/n}$ such that $F=(F_{1/n})^{*n}$, where ``$*$'' is the convolution, i.e. $F$ is $n$-fold convolution power of $F_{1/n}$.  It is well known that $F$ is infinitely divisible if and only if its characteristic function 
\begin{eqnarray*}
	f(t)\defeq\int\limits_{\R} e^{itx} \dd F(x),\quad t\in\R,
\end{eqnarray*}
admits the following representation
\begin{eqnarray}\label{repr_f}
	f(t) = \exp\Biggr\{it\gamma_\tau - \frac{\sigma^2 t^2}{2} + \int\limits_{\R\setminus\{0\}} \bigl(e^{itu}-1- \tfrac{it}{\tau}\sin (\tau u)\bigr)\dd \Lambda(u)\Biggr\}, \quad t\in\R,
\end{eqnarray} 
for any $\tau>0$ and some $\gamma_\tau\in\R$, $\sigma \geqslant 0$, and function $\Lambda$ (the \textit{L\'evy spectral function}), which is  non-decreasing on every interval $(-\infty,0)$ and $(0,+\infty)$, and it satisfies
\begin{eqnarray}\label{eq_Lambda_zerooninfty}
	\lim_{u \to -\infty}\Lambda(u) = \lim_{u\to+\infty}\Lambda(u)  = 0,
\end{eqnarray}
and also
\begin{eqnarray*}
	\int\limits_{0<|u|<\delta} x^2 \dd\Lambda(u)<+\infty,\quad\delta>0.
\end{eqnarray*}
Here $\gamma_\tau$ depend on $\tau$, but $\sigma$ and $\Lambda$ are not.  We use $u\mapsto \tfrac{1}{\tau}\sin (\tau u)$ as the ``centering function'' in the integral in \eqref{repr_f} following to Zolotarev \cite{Zolot1} and \cite{Zolot2}. This function can be replaced by any bounded Borel function $c:\R\to\R$, which satysfies
$c(0)=0$ and $(c(x)-x)/x^2\to 0$ as $x\to 0$ (see \cite{LindPanSato}).

It is well known that the properties of infinitely divisible $F$ (continuity, moments, and  support)  can be described in terms of its triplet $(\gamma_\tau, \sigma^2,\Lambda)$ (see \cite{Steutel}).

It turns out that there are probability laws, whose characteristic functions admit representation of type \eqref{repr_f} but with non-monotonic spectral function $\Lambda$.  The examples of such laws can be found in the classical monographs   \cite{GnedKolm}, \cite{LinOstr}, and  \cite{Lukacs}. The corresponding class  was introduced by Lindner and Sato \cite{LindSato}.  Following them, a distribution function $F$ is called \textit{quasi-infinitely divisible} if its characteristic function $f$ admits the representation \eqref{repr_f} with some (any) $\tau>0$ and some  $\gamma_\tau\in\R$, $\sigma \geqslant 0$, and function $\Lambda$, which has finite total variation on every interval $(-\infty,-r]$ and $[r,\infty)$, $r>0$, and it satisfies \eqref{eq_Lambda_zerooninfty} and
\begin{eqnarray*}
	\int\limits_{0<|u|<\delta} x^2 \dd|\Lambda|(u)<+\infty,\quad\delta>0.
\end{eqnarray*}
Such probability laws find interesting applications in theory of stochastic processes (see  \cite{LindSato} and \cite{Pass}), number theory (see \cite{Nakamura}),  physics (see \cite{ChhDemniMou} and \cite{Demni}), insurance mathematics (see  \cite{ZhangLiuLi}). 

The first detailed analysis of  quasi-infinitely divisible distribution functions on $\R$ was performed in \cite{LindPanSato} (the multivariate case is considered in the recent papers \cite{BergKutLind} and \cite{BergLind}).
In these works the authors studied  questions concerning  supports, moments, continuity, and weak convergence. The most complete results were obtained for probability laws on the set of integers $\Z$ (with extension for lattice laws).  The following important fact was stated in \cite{LindPanSato}.
\begin{Theorem}\label{th_QIZ}
	Let $F$ be the distribution function of probability law on $\Z$ with characteristic function $f$. Then $F$ is quasi-infinitely divisible if and only if $f$ does not have zeroes on the real line. In that case, $f$ admits the representation
	\begin{eqnarray}\label{th_QIZ_repr}
		f(t)=\exp\biggl\{it\gamma+  \sum_{k\in\Z\setminus\{0\}}\lambda_{k} (e^{itk}-1)\biggr\},\quad t\in\R,
	\end{eqnarray}
	where $\gamma\in\Z$, $\lambda_k\in\R$, $k\in\Z\setminus\{0\}$, and $\sum_{k\in\Z\setminus\{0\}}|\lambda_{k}|<\infty$.
\end{Theorem}
Here we omitted some detailes from the original formulation of this theorem. It is not difficult to check that \eqref{th_QIZ_repr} can be written in the form \eqref{repr_f} (we will do the similar derivation in the next section). 

In \cite{LindPanSato} the authors also obtained a criterion of weak convergence for sequences of quasi-infinitely divisible laws  to a fixed law from this class in terms of parameters of \eqref{th_QIZ_repr}. These results were complemented by compactness theorems in  \cite{Khartov}. It is known that weak and total variation convergences are equivalent for probability laws on $\Z$. Therefore these results can be considered  within the convergence in variation.

In this paper we generalize the above mentioned results. Namely, we state similar facts for rather wide class of discrete probability laws on $\R$, whose characteristic functions are separated from zero. In particular, this class includes arbitrary discrete infinitely divisible laws and lattice probability laws, whose characteristic functions have no zeroes on the real line. The first step in this direction was made in \cite{KhartAlexeev}.

\section{Representation theorem}
Let us consider an arbitrary discrete  distribution function
\begin{eqnarray*}
	F(x) \defeq \sum_{k=1}^{\infty}p_{x_k} I_{x_k}(x), \quad x \in \R,
\end{eqnarray*} 
where  $x_k$, $k\in\N$ (the set of positive integers), are distinct real numbers, $p_{x_k}\geqslant 0$, $k\in \N$, and $\sum_{k=1}^\infty p_{x_k}=1$.  Here and below $I_c(x)=1$ if $x\geqslant c$ and $I_c(x)=0$ if $x<c$ for any fixed $c\in\R$. Let $X$ be the set of all points of growth of $F$, i.e. $X\defeq \{x_k:  p_{x_k}>0,\, k\in\N \}\ne \varnothing$. So we can write
\begin{eqnarray*}
	F(x)=\sum_{x_k\in X} p_{x_k} I_{x_k}(x), \quad x \in \R.
\end{eqnarray*} 
Thus $F$ corresponds to a discrete probability law concentrated on the countable set $X$ of real numbers.

Let $f$ be the characteristic function of $F$. Then 
\begin{eqnarray*}
	f(t) = \sum_{k=1}^{\infty}p_{x_k} e^{it x_k}=\sum_{x_k\in X} p_{x_k} e^{it x_k}, \quad t \in \R.
\end{eqnarray*}
So $f$ is an almost periodic function with absolutely convergent Fourier series.

We focus on the case, when $f$ is separated from zero. Namely, we always assume that
$$
\text{\textit{there exists}\quad} \mu>0\quad\text{\textit{such that}} \quad|f(t)|\geqslant \mu>0\text{\quad \textit{for all}\quad} t\in\R. \eqno{(S)}
$$
Condition $(S)$ is equivalent to the relation $\inf_{t\in\R}|f(t)|>0$. 

Let $\D$ denote the class of all discrete distribution functions. Let $\DS$  be the subclass of $\D$,  which consists of distribution functions, whose characteristic functions satisfy condition $(S)$.  The class $\DS$ is rather wide. For example, if there exists $x_{k_0}$ such that $p_{x_{k_0}}>1/2$, then $F\in\DS$. This is seen from the following estimate:
\begin{eqnarray*}
	|f(t)|&=&\biggl|p_{x_{k_0}} e^{it x_{k_0}}+\sum_{k\in\N:\, k\ne k_0}p_{x_k} e^{it x_k}\biggr|
	\geqslant p_{x_{k_0}}-\sum_{k\in\N:\, k\ne k_0}p_{x_k}\\
	&&{}\qquad=p_{x_{k_0}}-(1-p_{x_{k_0}})=2p_{x_{k_0}}-1>0.
\end{eqnarray*}
Moreover, $\DS$ contains all  distribution functions of lattice probability laws, whose characteristic functions  have no zeroes on the real line, i.e. if
$X\subset\{a+b l: l\in\Z\}$ with some $a\in\R$ and $b>0$ and also $f(t)\ne 0$, $t\in\R$, then $F\in \DS$. Indeed, in this case $|f(\cdot)|$ is a continuous $2\pi/b$-periodic function and  its infimum is attained at some point in the period segment. Thus to have no zeroes on the real line means that $\inf_{t\in\R}|f(t)|>0$. 

It is important to note that all discrete infinitely divisible distribution functions belongs to $\DS$. Let us show this. It is known (see \cite{BlumRos}) that discrete distribution function $F$ is infinitely divisible if and only if its characteristic function
$f$ admits the following representation
 \begin{eqnarray}\label{eq_infdivf}
 	f(t) = \exp\biggr\{ it\gamma + \sum_{k=1}^\infty\lambda_{u_k}(e^{itu_k}-1)\biggr \}, \quad t\in\R,
\end{eqnarray}
where $\gamma\in\R$, $u_k\in\R$ and $\lambda_{u_k}\geqslant 0$ for all $k\in\N$, and $\sum_{k=1}^\infty\lambda_{u_k}<\infty$. Therefore we have the estimate
\begin{eqnarray*}
	|f(t)|=\exp\biggr\{ \sum_{k=1}^\infty\lambda_{u_k}(\cos(tu_k)-1)\biggr \}\geqslant \exp\biggl\{-2 \sum_{k=1}^\infty\lambda_{u_k}\biggr \}>0,\quad t\in\R.
\end{eqnarray*}
We see that $f$ satisfies $(S)$, i.e. $F\in \DS$.

There exist examples of characteristic functions of probability laws admitting representation \eqref{eq_infdivf}, where some of $\lambda_{u_k}$ can be negative within condition  $\sum_{k=1}^\infty|\lambda_{u_k}|<\infty$ (see \cite{GnedKolm} pp. 81--83,  \cite{LinOstr} p. 165, and \cite{Lukacs} p. 123--124). For such characteristic functions we similarly get the estimate
\begin{eqnarray*}
	|f(t)|\geqslant \exp\biggl\{-2 \sum_{k=1}^\infty|\lambda_{u_k}|\biggr \}>0,\quad t\in\R.
\end{eqnarray*}
Thus the corresponding distribution functions belong to $\DS$. In the result below we show that the converse fact  is also true, i.e. the characteristic function of every $F\in\DS$  admits representation of type  \eqref{eq_infdivf} with $\gamma\in\R$, $\lambda_{u_k}\in\R$, $k\in\N$, and $\sum_{k=1}^\infty|\lambda_{u_k}|<\infty$.

Howeever, it is clear that $\DS\ne\D$. For example, the distribution function $G(x)\defeq \tfrac{1}{2} I_0(x)+\tfrac{1}{2} I_1(x)$, $x\in\R$, does not belong to $\DS$. Indeed, its characteristic function $g(t)=\tfrac{1}{2}+\tfrac{1}{2} e^{it}$, $t\in\R$, has zeroes in the points $\pi (2k+1)$, $k\in\Z$. There are exist more interesting examples. Namely, let us consider the distribution function $H(x)\defeq\tfrac{1}{2} I_0 + p_\alpha I_{\alpha}(x)+ p_1 I_1(x)$, $x\in\R$, where $\alpha$ is an irrational number from $(0,1)$, the weights $p_\alpha$ and $p_1$ are strictly positive, and  $p_\alpha+p_1=1/2$. Its characteristic function $h(t)= \tfrac{1}{2}+p_\alpha e^{it\alpha}+ p_1 e^{it}$, $t\in\R$, does not have zeroes, and it is not separated from zero, i.e. $h(t)\ne 0$, $t\in\R$, and $\inf_{t\in\R}|h(t)|=0$. Thus we have $H\notin \DS$.  This and other similar particular results  were obtained in \cite{KhartAlexeev}.

For sharp formulation of the result below we  need to introduce the set of all finite $\Z$--linear combinations of elements from a set $Y\subset \R$:
\begin{eqnarray*}
	\la Y\ra\defeq \biggl\{\sum_{k=1}^m c_k y_k:\,   c_k \in \Z,\, y_k\in Y,\, m\in\N\biggr\}.
\end{eqnarray*}
In other words, $\la Y\ra$ is a module over the ring $\Z$ with the generating set $Y$. It easily seen that, in particular, $Y\subset\la Y\ra$, $0\in \la Y\ra$, and $k y\in\la Y \ra$ for any  $k\in\Z$, $y\in Y$. If $Y\ne\{0\}$, then $\la Y\ra$ is an infinite countable set. 
  
We always set below that sums over the null set equal zero. 
  
\begin{Theorem}\label{th_Repr}
	Under the assumption $(S)$, the characteristic function $f$ admits the following representation
 	\begin{eqnarray}\label{th_Repr_eq}
 		f(t) = \exp\biggr\{ it\gamma + \sum_{u \in \la X\ra\setminus \{0\}}\lambda_{u}(e^{itu}-1)\biggr \}, \quad t\in\R,
 	\end{eqnarray}
 	where $\gamma \in \la X\ra$, $\lambda_u \in \R$ for all $u\in\la X\ra\setminus \{0\}$, and $\sum_{u \in \la X\ra\setminus \{0\}}|\lambda_u| < \infty$.
\end{Theorem}

Before the proof of Theorem \ref{th_Repr} we will do some remarks.

According to this theorem, if $F\in\DS$, i.e. $f$ satisfies $(S)$, then $f$ admits representation \eqref{th_Repr_eq}. Let $\Arg f(t)$ be the argument of $f(t)$  (i.e.  $f(t)=\exp\{i\Arg f(t)\}$) that is uniquely defined for $t\in\R$ by continuity with the condition $\Arg f(0)=0$.  Then
\begin{eqnarray*}
	\Arg f(t)=t\gamma+\sum_{u \in \la X\ra\setminus \{0\}}\lambda_{u}\sin(tu),\quad t\in\R.
\end{eqnarray*}
Due to $\sum_{u \in \la X\ra\setminus \{0\}}|\lambda_u| < \infty$, we have
\begin{eqnarray}\label{conc_gamma}
	\gamma=\lim\limits_{T\to\infty} \dfrac{\Arg f(T)}{T}=\lim\limits_{T\to\infty} \dfrac{\Arg f(T)-\Arg f(0)}{T-0}. 
\end{eqnarray}
Thus $\gamma$ is a ``mean motion'' of the argument of $f$ on $[0,\infty)$. 

Let us fix arbitrarily $\tau>0$ and introduce  ``mean motion'' on the interval $[0,\tau]$:
\begin{eqnarray}\label{def_gammatau}
	\gamma_\tau\defeq \dfrac{\Arg f(\tau)}{\tau}=\gamma+\tfrac{1}{\tau}\!\!\!\sum_{u \in \la X\ra\setminus \{0\}}\lambda_{u}\sin(\tau u).
\end{eqnarray}
Then we write
\begin{eqnarray*}
	f(t)=\exp\biggl\{it\gamma_\tau+ \sum_{u \in \la X\ra\setminus \{0\}}\lambda_{u} (e^{itu}-1-i\tfrac{t}{\tau}\sin (\tau u)) \biggr\}, \quad t\in\R.
\end{eqnarray*}
Let us define
\begin{eqnarray*}
	\Lambda(u)=\sum\limits_{k\in\N: u_k\leqslant u} \lambda_{u_k},\quad u<0,\quad\text{and}\quad \Lambda(u)=-\sum\limits_{k\in\N: u_k> u} \lambda_{u_k},\quad u>0.
\end{eqnarray*}
Hence we have
\begin{eqnarray*}
	f(t)=\exp\Biggr\{it\gamma_\tau  + \int\limits_{\R\setminus\{0\}} \bigl(e^{itu}-1-\tfrac{it}{\tau} \sin (\tau u)\bigr)\dd \Lambda(u)\Biggr\}, \quad t\in\R.
\end{eqnarray*}
So we proved the following fact.
\begin{Remark}\label{rem_th_Repr}
	If $F\in\DS$, then it is quasi-infinitely divisible with the triplet $(\gamma_\tau, 0, \Lambda)$.
\end{Remark}
It is easily seen that the function $\Lambda$ satisfies all the required conditions.\\

Note that the centering by $u\mapsto \tfrac{1}{\tau} \sin (\tau u)$ in the integral has the advantage that the shift parameter $\gamma_\tau$ of the triplet gets additional meaning according to formulas \eqref{conc_gamma} and \eqref{def_gammatau}.

For fixed $\tau>0$ the parameter $\gamma_\tau$ and the spectral function $\Lambda$ are uniquely determined from $f$ (see \cite{GnedKolm}, p. 80). Therefore  $\gamma$ and the sequence  $\lambda_{u}$, $u\in \langle X\rangle \setminus \{0\}$, are also uniquely determined in \eqref{th_Repr_eq}.

According to Remark \ref{rem_th_Repr} and comments above concerning lattice probability laws, Theorem \ref{th_Repr} is a generalization of Theorem \ref{th_QIZ}.

It is very interesting to know whether or not every discrete quasi-infinitely divisible distribution function belongs to $\DS$. It is known that characteristic functions of quasi-infinitely divisible distribution functions have no zeroes on the real line (see \cite{LindPanSato}). So the problem is reduced to the following question. Are there exist discrete quasi-infinitely divisible distribution functions with characteristic functions that have no zeroes on the real line and that are not separated from zero? There are only particular results. For example, the above mentioned distribution function $H$ is not quasi-infinitely divisible that was shown in \cite{KhartAlexeev}.\\

\textbf{Proof of Theorem \ref{th_Repr}.} The proof is devided into two steps: for the case of finite $X$ and for the general case. 

1. Suppose that $X$ is finite. Without loss of generality we can assume that $X=\{x_1,\ldots, x_n\}$.  Hence $f(t)=\sum_{k=1}^n p_{x_k} e^{it x_k}$, $t \in \R$.  The case $n=1$ and $x_1=0$ is trivial: $\la X\ra=\{0\}$, $\gamma=0$. We next exclude this case. Let $\alpha_1, \ldots, \alpha_d\in X$ be a basis in $X\ne\{0\}$ over  $\Z$, i.e. for every $k=1,\ldots, n$ there is the unique representation  $x_k=\sum_{j=1}^d c_{k,j} \alpha_j$ with $c_{k,j}\in\Z$. Here the numbers $\alpha_1, \ldots, \alpha_d$ are linearly independent over $\Z$, that is the equation $l_1\alpha_1+\ldots +l_1\alpha_d=0$ with $l_1,\ldots, l_d\in\Z$ holds if and only if $l_1=l_2=\ldots=l_d=0$. It is easy to check that
\begin{eqnarray}\label{eq_moduleX_alpha}
	\la X\ra= \biggl\{\sum_{j=1}^d l_j \alpha_j:\,   l_j \in \Z\biggr\}.
\end{eqnarray}

Using the representations for $x_k$,  we consider the function
\begin{eqnarray*}
	\varphi(t_1,\ldots, t_d)=\sum_{k=1}^n p_{x_k} \exp\biggl\{i\sum\limits_{j=1}^d c_{k,j} \alpha_j t_j\biggr\},\quad t_1,\ldots, t_d \in \R. 
\end{eqnarray*}
The function $\varphi$ is a continuous periodic function  with periods $2\pi/\alpha_1,\ldots, 2\pi/\alpha_d$ over variables $t_1,\ldots,t_d$, correspondingly.  The function $f$ is the diagonal function for $\varphi$:
\begin{eqnarray*}
	\varphi(t,\ldots,t)=\sum_{k=1}^n p_{x_k} \exp\biggl\{it\sum\limits_{j=1}^d c_{k,j} \alpha_j \biggr\}=\sum_{k=1}^n p_{x_k} e^{it x_k}=f(t),\quad t\in\R.
\end{eqnarray*}
Therefore the set of values $\{f(t): t\in\R\}$ is everywhere dense in the set $\{\varphi(t_1, \ldots, t_d): t_1,\ldots, t_d\in \R\}$ (see \cite{LevitanZhikov}, Theorem 5,  p. 46). Hence, by the assumption $(S)$, we have
\begin{eqnarray}\label{ineq_phi_mu}
	|\varphi(t_1, \ldots, t_d)| \geqslant \mu >0, \quad t_1, \ldots, t_d\in\R.
\end{eqnarray} 
We now consider  the following function
\begin{eqnarray*}
	\tilde\varphi(t_1,\ldots,t_d)\defeq \varphi(t_1/\alpha_1,\ldots, t_d/\alpha_d)=\sum_{k=1}^n p_{x_k} \exp\biggl\{i\sum\limits_{j=1}^d c_{k,j}  t_j\biggr\},\quad t_1,\ldots, t_d \in\R.
\end{eqnarray*}
It is  the characteristic function for the probabilty law on the vectors $(c_{k,1},\ldots, c_{k,d})\in\Z^d$ with weights $p_{x_k}$, $k=1,\ldots, n$. Due to \eqref{ineq_phi_mu}, we have  $|\tilde\varphi(t_1, \ldots, t_d)| \geqslant \mu >0$,  $t_1, \ldots, t_d\in\R$. By Theorem 3.2 in \cite{BergLind}, the function $\tilde{\varphi}$ admits the following representation for any $t_1, \ldots, t_d \in \R$:
\begin{eqnarray*}
	\tilde\varphi(t_1, \ldots, t_d) &=&\exp\biggl\{ i(t_1\gamma_1  + \ldots+ t_d \gamma_d)\\
	 &&{}\qquad+\!\!\! \sum_{(l_1, \ldots, l_d) \in \Z^d\setminus\{0\}}\!\!\!\lambda_{l_1, \ldots, l_d} \Bigr(\exp\bigl\{i(t_1l_1+\ldots+ t_d l_d )\bigr\}-1\Bigr)\biggr\},
\end{eqnarray*}
where $\gamma_1, \ldots, \gamma_d \in \Z$, $\lambda_{l_1, \ldots, l_d} \in \R$ for  $(l_1, \ldots, l_d)\in \Z^d\setminus \{0\}$, and $\sum_{(l_1, \ldots, l_d) \in \Z^d\setminus\{0\}}|\lambda_{l_1, \ldots, l_d}| < \infty$. 

Next, observe that
\begin{eqnarray*}
 \tilde\varphi(t\alpha_1 , \ldots, t\alpha_d )=\varphi(t, \ldots,  t)=f(t),\quad t\in\R.
\end{eqnarray*}
Then for $t\in\R$ we have
\begin{eqnarray*}
	f(t) =\exp\biggl\{ it(\gamma_1\alpha_1  + \ldots+ \gamma_d\alpha_d ) +\!\!\! \sum_{(l_1, \ldots, l_d) \in \Z^d\setminus\{0\}}\!\!\!\lambda_{l_1, \ldots, l_d} \Bigr(\exp\bigl\{it(l_1\alpha_1+\ldots+ l_d \alpha_d)\bigr\}-1\Bigr)\biggr\}.
\end{eqnarray*}
Due to \eqref{eq_moduleX_alpha}, this representation we can write in the form \eqref{th_Repr_eq}, where $\gamma\defeq \gamma_1\alpha_1  + \ldots+ \gamma_d\alpha_d\in\la X\ra$, and $\lambda_u\defeq \lambda_{l_1, \ldots, l_d}$ for $u=l_1\alpha_1+\ldots+ l_d \alpha_d\in\la X\ra\setminus\{0\}$ with $(l_1, \ldots, l_d) \in \Z^d\setminus\{0\}$.\\

2. We now turn to the general case for $X$. So $f(t)=\sum_{k=1}^{\infty} p_{x_k}e^{itx_k}$, $t\in\R$. Since $X\ne \varnothing$, without loss of generality we can suppose that $x_1\in X$, i.e. $p_{x_1}>0$. We will uniformly approximate $f$ by the following characteristic functions:
\begin{eqnarray*}
	f_n(t)\defeq \sum_{k=1}^{n}  p_{n,x_k} e^{itx_k},\quad  t\in\R,\quad n\in\N,
\end{eqnarray*}
where
\begin{eqnarray*}
	p_{n,x_k}\defeq \dfrac{p_{x_k}}{\sum_{m=1}^{n} p_{x_m}},\quad k=1,\ldots, n,\quad n\in\N.
\end{eqnarray*}
Here $\sum_{m=1}^{n} p_{x_m}\geqslant p_{x_1}>0$, $n\in\N$.   Let us estimate the approximation error for every $n\in\N$:
\begin{eqnarray*}
	\sup\limits_{t\in\R}|f(t)-f_n(t)|&=&\sup\limits_{t\in\R}\biggl| \sum_{k=1}^n (p_{x_k}-p_{n,x_k})e^{itx_k} +\sum_{k=n+1}^\infty p_{x_k}e^{itx_k}\biggr|\\
	&\leqslant& \sum_{k=1}^n |p_{x_k}-p_{n,x_k}|+\sum_{k=n+1}^\infty p_{x_k},
\end{eqnarray*}
where, due to the equality $\sum_{m=1}^{\infty} p_{x_m}=1$, we have
\begin{eqnarray*}
 \sum_{k=1}^n |p_{x_k}-p_{n,x_k}|= \sum_{k=1}^n p_{x_k}\biggl(\dfrac{1}{\sum_{m=1}^{n} p_{x_m}}-1\biggr)
	= \sum_{k=1}^n p_{x_k}\biggl(\dfrac{\sum_{m=n+1}^{\infty} p_{x_m}}{\sum_{m=1}^{n} p_{x_m}}\biggr)=\sum_{m=n+1}^{\infty} p_{x_m}.
\end{eqnarray*}
Thus we obtain
\begin{eqnarray*}
	\sup\limits_{t\in\R}|f(t)-f_n(t)|\leqslant 2\sum_{k=n+1}^\infty p_{x_k},\quad n\in\N.
\end{eqnarray*}
From this we conclude that 
\begin{eqnarray*}
	\sup_{t\in\R}|f(t)-f_n(t)|\to 0,\quad n\to\infty.
\end{eqnarray*}
Let us take $\mu>0$ from $(S)$, and fix arbitrarily $\e\in(0,1/4)$. Then there exists $n_\e\in\N$ such that
\begin{eqnarray}\label{conc_supinft_ffn}
	\sup_{t\in\R} |f(t)-f_n(t)|\leqslant\e \mu,\quad n\geqslant n_\e.
\end{eqnarray} 
From this we see that $f_n(t)\ne 0$, $t\in\R$, for every $n\geqslant n_\e$. Indeed,
\begin{eqnarray*}
	|f_n(t)|\geqslant |f(t)|-|f(t)-f_n(t)|\geqslant \inf_{t\in\R} |f(t)| - \sup_{t\in\R} |f(t) -f_n(t)|,
\end{eqnarray*}
and, due to $(S)$ and \eqref{conc_supinft_ffn}, we have
\begin{eqnarray}\label{ineq_fn_e_mu}
	|f_n(t)|\geqslant(1-\e)\mu >0,\quad t\in\R.
\end{eqnarray} 

Choose $n\geqslant n_\e$ and represent $f(t)=f_n(t)\cdot R_n(t)$ with   $R_n(t)\defeq f(t)/f_n(t)$, $t\in\R$. Since $f$, $f_n$, $R_n$ are continuous functions without zeroes  and they  equal $1$ at $t=0$, we can proceed to  \textit{the distinguished logarithms}: 
\begin{eqnarray}\label{eq_Lnf}
	\Ln f(t)=\Ln f_n(t) + \Ln R_n(t),\quad t\in\R,
\end{eqnarray}
which are uniquely defined by continuity with the conditions that they equal $0$ at $t=0$ (see \cite{Khinch} p. 14--16  or \cite{Sato} p. 2). We next consider $\Ln f_n$ and  $\Ln R_n$ separately.

By the part 1., for the function $f_n$ we have
\begin{eqnarray*}
	f_n(t)=\exp\biggl\{it\gamma_n+
	\sum_{u\in\la X_n\ra\setminus\{0\}}\lambda_{n,u} (e^{itu}-1)\biggr\},\quad t\in\R,
\end{eqnarray*}
where $X_n=\bigl\{x_k: p_{x_k}>0,\, k=1,\ldots,n \bigr\}$, $\gamma_n\in\la X_n\ra$, $\lambda_{n,u}\in\R$ for all $u\in\la X_n \ra\setminus\{0\}$, and $\sum_{u\in\la X_n\ra\setminus\{0\}} |\lambda_{n,u}|<\infty$. Since the function in the exponent is continuous and it equals $0$ at $t=0$, we have
\begin{eqnarray*}
	\Ln f_n(t)=it\gamma_n+
	\sum_{u\in\la X_n\ra\setminus\{0\}}\lambda_{n,u} (e^{itu}-1),\quad t\in\R.
\end{eqnarray*}
Setting $\lambda_{n,0}\defeq -\sum_{u\in\la X_n\ra\setminus\{0\}} \lambda_{n,u}$ , we represent $f_n$ in the following form
\begin{eqnarray*}
	\Ln f_n(t)=it\gamma_n+
	\sum_{u\in\la X_n\ra}\lambda_{n,u} e^{itu},\quad t\in\R.
\end{eqnarray*}
Observe that $X_n\subset X$, and hence $\la X_n\ra\subset \la X\ra$. Then we can write
\begin{eqnarray}\label{eq_fn_moduleX}
	\Ln f_n(t)=it\gamma_n+
	\sum_{u\in\la X\ra}\lambda_{n,u} e^{itu},\quad t\in\R,
\end{eqnarray}
where for every $u\in \la X\ra\setminus \la X_n\ra$ we define $\lambda_{n,u}\defeq 0$  if $\la X\ra\setminus \la X_n\ra\ne \varnothing$.

Next, we consider the function $\Ln R_n$. Observe that
\begin{eqnarray}\label{eq_Lnffn}
	\Ln R_n(t)=\ln\biggl(1+\dfrac{f(t)-f_n(t)}{f_n(t)}\biggr),\quad t\in\R,
\end{eqnarray}
(the latter is the principal value of the logarithm) because due to \eqref{conc_supinft_ffn} and \eqref{ineq_fn_e_mu} 
\begin{eqnarray}\label{conc_supt_ffn}
	\sup_{t\in\R} \biggl|\dfrac{f(t) -f_n(t)}{f_n(t)}\biggr| \leqslant \dfrac{\sup_{t\in\R}|f(t) -f_n(t)|}{\inf_{t\in\R}|f_n(t)|}\leqslant \dfrac{\e}{1-\e}<1,
\end{eqnarray}
and the function in the right-hand side of \eqref{eq_Lnffn} is continuous and equals $0$ at $t=0$.
Using this, we get the decomposition
\begin{eqnarray*}
 \Ln R_n(t)=\sum\limits_{m=1}^{\infty} \dfrac{(-1)^{m-1}}{m} \biggl(\dfrac{f(t)-f_n(t)}{f_n(t)}\biggr)^m,\quad t\in\R.
\end{eqnarray*}
It yields the estimate
\begin{eqnarray*}
	\sup_{t\in\R}|\Ln R_n(t)|\leqslant \sum\limits_{m=1}^{\infty} \dfrac{1}{m}\, \sup_{t\in\R}\biggl|\dfrac{f(t)-f_n(t)}{f_n(t)}\biggr|^m\leqslant \sum\limits_{m=1}^{\infty} \dfrac{1}{m} \biggl(\dfrac{\e}{1-\e}\biggr)^m.
\end{eqnarray*}
Since $\e\in(0,1/4)$, here we have
\begin{eqnarray*}
	\sum\limits_{m=1}^{\infty} \dfrac{1}{m} \biggl(\dfrac{\e}{1-\e}\biggr)^m\leqslant\sum\limits_{m=1}^{\infty}  \biggl(\dfrac{\e}{1-\e}\biggr)^m=\dfrac{\tfrac{\e}{1-\e}}{1-\tfrac{\e}{1-\e}}=\dfrac{\e}{1-2\e}<2\e.
\end{eqnarray*}
Thus we obtain
\begin{eqnarray}\label{ineq_supRn}
	\sup_{t\in\R}|\Ln R_n(t)|<2\e.
\end{eqnarray}

Let us consider the function $(f-f_n)/f_n$ from \eqref{eq_Lnffn}. It is clear that $f-f_n$ is an almost periodic function with absolutely convergent Fourier series. Due to \eqref{ineq_fn_e_mu}, the function $1/f_n$ is also almost periodic (see \cite{Bohr} p. 35 or \cite{Levitan} p. 21),  and it has absolutely convergent Fourier series (see \cite{Gelfand} p. 175). These facts can be also shown using \eqref{eq_fn_moduleX}. Therefore the function $(f-f_n)/f_n$ is almost periodic with absolutely convergent Fourier series.

Let us return to \eqref{eq_Lnffn}. Since the function $z\mapsto\ln(1+z)$, $z\in\CC$ (the set of complex numbers), is analytic on the unit disk, due to \eqref{conc_supt_ffn} and the above remark, we get that $\Ln R_n$ is almost periodic with absolute convergent Fourier series (see \cite{Gelfand} p. 175):
\begin{eqnarray}\label{eq_LnRn}
	\Ln R_n(t)=\sum\limits_{u\in \Delta_n} \beta_{n,u} e^{it u},\quad t\in\R,
\end{eqnarray}
where $\Delta_n$ is a countable set of real numbers (and, probably, $\Delta_n\ne \la X\ra$), $\beta_{n,u}\in \CC$ for $u\in \Delta_n$, and $\sum_{u\in \Delta_n}|\beta_{n,u}|<\infty$. Without loss of generality we can assume that either $\Delta_n=\varnothing$ (the case $f=f_n$) and the sum equals zero in \eqref{eq_LnRn}, or  $\Delta_n\ne\varnothing$ and $\beta_{n,u}\ne 0$, $u\in\Delta_n$.

We now return to the function $\Ln f$ and \eqref{eq_Lnf}.  The formulas \eqref{eq_fn_moduleX} and \eqref{eq_LnRn} yield 
\begin{eqnarray*}
	\Ln f(t)=it\gamma_n +\sum_{u\in\la X\ra}\lambda_{n,u} e^{itu}+ \sum\limits_{u\in \Delta_n} \beta_{n,u} e^{it u},\quad t\in\R.
\end{eqnarray*}
Recall that this is valid for every $n\geqslant n_\e$. Since $\sum_{u\in\la X\ra}|\lambda_{n,u}|<\infty$ and $\sum_{u\in \Delta_n} |\beta_{n,u}|<\infty$, $n\geqslant n_\e$, it is true that
\begin{eqnarray*}
	\lim\limits_{T\to\infty} \dfrac{\Ln f(T)}{iT}= \gamma_n,\quad n\geqslant n_\e.
\end{eqnarray*}
Consequently, $\gamma_n$ are equal for $n\geqslant n_\e$, and  we denote their value by $\gamma$. Since $\gamma_n\in\la X_n\ra \subset \la X\ra$, we have $\gamma \in \la X \ra$. 

We next consider the function
\begin{eqnarray}\label{eq_Lnf_gamma_series}
	\Ln f(t)-it\gamma=\sum_{u\in\la X\ra}\lambda_{n,u} e^{itu}+ \sum\limits_{u\in \Delta_n} \beta_{n,u} e^{it u},\quad t\in\R,\quad n\geqslant n_\e.
\end{eqnarray}
It is an almost periodic function with absolutely convergent Fourier series. We define $\Delta=\cup_{n=n_\e}^\infty\Delta_n$, and we introduce the Fourier coefficients
\begin{eqnarray}\label{def_lambdav}
	\lambda_v\defeq \lim\limits_{T\to\infty} \dfrac{1}{2T} \int\limits_{-T}^{T} \bigl(\Ln f(t)-it\gamma  \bigr)e^{-it v}\dd t,\quad v\in \la X\ra \cup \Delta.
\end{eqnarray}
According to \eqref{eq_LnRn} and \eqref{eq_Lnf_gamma_series},  for any $n\geqslant n_\e$ we have
\begin{eqnarray*}
	\lambda_v=\lim\limits_{T\to\infty} \dfrac{1}{2T} \int\limits_{-T}^{T} \biggl(\,\sum_{u\in\la X\ra}\lambda_{n,u} e^{itu} \biggr)e^{-it v}\dd t+ \lim\limits_{T\to\infty} \dfrac{1}{2T} \int\limits_{-T}^{T} \Ln R_n(t)e^{-it v}\dd t,\quad v\in \la X\ra \cup \Delta.
\end{eqnarray*}
Due to \eqref{ineq_supRn}, for any $n\geqslant n_\e$ 
\begin{eqnarray*}
	\Biggl|  \dfrac{1}{2T} \int\limits_{-T}^{T} \Ln R_n(t)e^{-it v}\dd t \Biggr|\leqslant \dfrac{1}{2T} \int\limits_{-T}^{T} | \Ln R_n(t)|\dd t\leqslant \sup_{t\in\R} | \Ln R_n(t)|\leqslant 2\e.
\end{eqnarray*}
Since $\e>0$ can be choosen arbitrarily small, we obtain
\begin{eqnarray*}
	\lambda_v=\lim_{n\to\infty}\lim\limits_{T\to\infty} \dfrac{1}{2T} \int\limits_{-T}^{T} \biggl(\,\sum_{u\in\la X\ra}\lambda_{n,u} e^{itu} \biggr)e^{-it v}\dd t,\quad v\in \la X\ra \cup \Delta.
\end{eqnarray*}
The series in the integral is absolutely convergent. Hence
\begin{eqnarray}\label{eq_lambdav}
	\lambda_v=\lim_{n\to\infty}\sum_{u\in\la X\ra}\lambda_{n,u}\biggl(\,\lim\limits_{T\to\infty} \dfrac{1}{2T} \int\limits_{-T}^{T}  e^{itu} e^{-it v}\dd t\biggr),\quad v\in \la X\ra \cup \Delta,
\end{eqnarray}
where
\begin{eqnarray}\label{eq_intT}
	\lim\limits_{T\to\infty} \dfrac{1}{2T} \int\limits_{-T}^{T}  e^{itu} e^{-it v}\dd t=\begin{cases}
		1, & u=v,\\
		0,& u\ne v.
	\end{cases}
\end{eqnarray}

Suppose that $\Delta\setminus \la X\ra \ne \varnothing$,  and $v\in\Delta\setminus \la X\ra$. Then, on the one hand,  from \eqref{eq_lambdav} and \eqref{eq_intT} we have $\lambda_v =0$. On the other hand, $v\in \Delta_n\setminus \la X\ra$ for some $n$. So from \eqref{eq_Lnf_gamma_series},  \eqref{def_lambdav}, and \eqref{eq_intT} we have $\lambda_v=\beta_{n,v}$ for such $n$, where $\beta_{n,v}\ne 0$ according to the convention above. Thus we obtain a contradiction.  Hence $\Delta\subset \la X\ra$. Therefore
\begin{eqnarray*}
	\Ln f(t)-it\gamma=\sum_{u\in\la X\ra}(\lambda_{n,u}+\beta_{n,u}) e^{itu},\quad t\in\R,\quad n\geqslant n_\e,
\end{eqnarray*}
where $\beta_{n,u}\defeq 0$, $u\in \la X\ra\setminus \Delta_n$. So $\lambda_{u}=\lambda_{n,u}+\beta_{n,u}$, $u\in\la X\ra$, $n\geqslant n_\e$, and, from \eqref{eq_lambdav} and \eqref{eq_intT}, we have $\lambda_u=\lim\limits_{n\to\infty}\lambda_{n,u}$, $u\in\la X\ra$. Since $\lambda_{n,u}\in \R$, we conclude that $\lambda_u\in\R$. Thus 
\begin{eqnarray*}
	\Ln f(t)-it\gamma=\sum_{u\in\la X\ra}\lambda_{u} e^{itu},\quad t\in\R,
\end{eqnarray*}
with $\lambda_u\in\R$ for $u\in\la X\ra$, and $\sum_{u\in\la X\ra}|\lambda_{u}|<\infty$.
Since $\bigl(\Ln f(t)-it\gamma\bigr)\bigr|_{t=0}=0$, we get $\lambda_0=-\sum_{u\in\la X\ra\setminus\{0\}}\lambda_{u}$. From this we obtain \eqref{th_Repr_eq}.  \quad $\Box$

\section{Limit and compactness theorems}

In this section we will apply the representations \eqref{th_Repr_eq} to obtain 
limit and compactness theorems  with convergence in variation to distribution functions  from $\DS$.

Let $(F_n)_{n\in\N}$ be a sequence of distribution functions and let $F$ be a distribution function. Recall that, by definition, the sequence $(F_n)_{n\in\N}$  \textit{converges in variation} to $F$ if $\|F_n -F\|\to 0$, $n\to\infty$, where $\|\cdot\|$ denotes the total variation on $\R$. We denote this as $F_n\tovar F$.
Here we are interested in the case, when  $F\in\DS$. So it is naturally to consider the following notions of compactness with class $\DS$.  We say that $(F_n)_{n\in\N}$ is $\DS$\textit{-relatively compact in variation} if every its subsequence contains a further subsequence  converging in variation to a distribution function from $\DS$. We say that $(F_n)_{n\in\N}$ is $\DS$\textit{-stochastically compact in variation} if every its subsequence contains a further subsequence  converging in variation to a nondegenerate distribution function from $\DS$. In particular, if a sequence of distribution functions converges in variation to some distribution function from  $\DS$, then the sequence is $\DS$-relatively compact in variation. If we know that the limit is nondegenerate, then the sequence is  $\DS$-stochastically compact in variation. The introduced notions of compactness are modified versions of well known relative (weak)  compactness (see \cite{Bill})  and stochastic compactness (introduced by Feller in \cite{Feller}), which are based on weak convergence.

It should be noted that the class $\DS$ is not closed under convergence in variation. Indeed, it is easily seen that the sequence 
\begin{eqnarray*}
	G_n(x)\defeq \bigl(\tfrac{1}{2}+\tfrac{1}{2+n}\bigr)\,I_0(x)+\bigl(\tfrac{1}{2}-\tfrac{1}{2+n}\bigr)\,I_1(x),\quad x\in\R,\quad n\in\N,
\end{eqnarray*}
 converges in variation to the law
\begin{eqnarray*}
	G(x)\defeq\tfrac{1}{2}\,I_0(x)+\tfrac{1}{2}\,I_1(x),\quad x\in\R.
\end{eqnarray*}
However, $G_n\in\DS$ for every $n\in\N$, but $G\notin\DS$ (see the examples from  Section 2). Here we also see that $(G_n)_{n\in\N}$ is not $\DS$-relatively (and hence $\DS$-stochastically) compact in variation.

We begin with the following observation.
\begin{Theorem}\label{th_Fn_D_F_DS}
	Suppose that $F_n\in\D$, $n\in\N$. If $(F_n)_{n\in\N}$ is $\DS$-relatively compact in variation, then  $F_n\in\DS$ for all sufficiently large $n\in\N$. 
\end{Theorem}

Due to this theorem, the conditions of $\DS$-relative compactness of $F_n\in \D$, $n\in\N$, can be formulated in terms of parameters from existing representations of type \eqref{th_Repr_eq} for characteristic functions of $F_n$ with all sufficiently large  $n\in\N$.  In particular, the theorem and this remark are valid for the sequences from $\D$ converging in variation to some $F\in\DS$.  

We fix $F\in \DS$. Let $f$ be its characteristic function with representation \eqref{th_Repr_eq}. Without loss of generality we suppose that the whole sequence $(F_n)_{n\in\N}$ is  from $\DS$. Let  $X_n$ be the  set of all points of growth of $F_n$ for every  $n\in\N$. Let $f_n$ be the characteristic function of $F_n$ for every $n\in\N$. According to Theorem \ref{th_Repr}, the functions $f_n$, $n\in\N$, admit representation of type \eqref{th_Repr_eq}. We set 
\begin{eqnarray*}
	f_n(t)=\exp\biggl\{it\gamma_n+  \sum_{u\in\la X_n\ra\setminus\{0\}}\lambda_{n,u} (e^{itu}-1)\biggr\},\quad t\in\R,\quad n\in\N, 
\end{eqnarray*}
where  $\gamma_n\in\la X_n\ra$, $\lambda_{n,u}\in\R$ for all $u\in\la X_n\ra\setminus\{0\}$, and $\sum_{u\in\la X_n\ra\setminus\{0\}}|\lambda_{n,u}|<\infty$.

We adopt the convention:
\begin{eqnarray}\label{conv_lambdau_lambdanu}
	\lambda_{u}\defeq 0,\quad u\notin \la X\ra,\quad\text{and}\quad\lambda_{n,u}\defeq 0,\quad u\notin \la X_n\ra,\quad n\in\N.
\end{eqnarray}

We propose the following criterion of convergence in variation.
\begin{Theorem}\label{th_Conv}
	 $F_n\tovar F$, $n\to\infty$, if and only if the following couple of conditions holds:
	\begin{itemize}
		\item[$(i)$] there exists $n_0\in\N$  such that  $\gamma_n=\gamma$, $n\geqslant n_0$;
		\item[$(ii)$] $\sum\limits_{\substack{u\in\la X_n \ra\cup \la X\ra,\\ u\ne0}}|\lambda_{n,u}-\lambda_{u}|\to 0$,  $n\to\infty$.
	\end{itemize}	 
\end{Theorem}
This criterion generalizes the corresponding result from \cite{LindPanSato} for probability laws on $\Z$. 

We next turn to compactness theorems. We need to introduce the following set
\begin{eqnarray*}
	\UnionXn=\bigcup_{n=1}^\infty \la X_n\ra=\{u_0,u_1,u_2,\ldots\}.
\end{eqnarray*}
Since always $0\in \la X_n\ra$, $n\in\N$, we set $u_0\defeq 0$.
\begin{Theorem}\label{th_RelComp}
	$(F_n)_{n\in\N}$ is $\DS$-relatively compact in variation if and only if the following ensemble of conditions holds:
	\begin{itemize}
		\item[$(i)$] $\gamma_n$, $n\in\N$, have finite number of values from $\UnionXn$;
		\item[$(ii)$] $\sup\limits_{n\in\N}\sum\limits_{u\in\la X_n \ra\setminus\{0\}}|\lambda_{n,u}|<\infty$;
		\item[$(iii)$] $\lim\limits_{N\to\infty}\sup\limits_{n\in\N}\sum\limits_{\substack{u_k\in\, \UnionXn\setminus\{0\}:\\ k>N}}|\lambda_{n,u_k}|=0$.
	\end{itemize}	 
\end{Theorem}

\begin{Theorem}\label{th_StochComp}
	$(F_n)_{n\in\N}$ is $\DS$-stochastically compact in variation if and only if it is $\DS$-relatively compact in variation, and 
	\begin{eqnarray}\label{th_StochComp_cond}
		\varliminf\limits_{n\to\infty}\sum\limits_{u\in\la X_n \ra\setminus\{0\}}|\lambda_{n,u}|>0.
	\end{eqnarray}	 
\end{Theorem}

Theorems \ref{th_RelComp} and \ref{th_StochComp} generalize the corresponding results from \cite{Khartov} concerning compactness  for probability laws on $\Z$.

We now turn to the proofs of proposed theorems.\\

\noindent
\textbf{Proof of Theorem \ref{th_Fn_D_F_DS}.} Suppose, contrary to our claim, that there exists a strictly increasing sequence  of positive integers $(n_j)_{j\in\N}$ such that every $F_{n_j}$ does not belong to $\DS$. Let us choose from $(F_{n_j})_{j\in\N}$ a subsequence $(F_{n_j'})_{j\in\N}$ that converges in variation to some distribution function  $F$ from $\DS$. Let $f_{n_j'}$ be the characteristic function of $F_{n_j'}$ for every $j\in\N$, and $f$ be the characteristic function of $F$. Then there exists $\mu>0$ such that $|f(t)|\geqslant \mu>0$ for all $t\in\R$.  Due to the estimate
\begin{eqnarray}\label{ineq_sup_absnorm}
	\sup_{t\in\R}|f_{n_j'}(t)-f(t)| =\sup_{t\in\R}\,\biggl| \int\limits_{\R}e^{itx} \dd (F_{n_j'}(x)-F(x))\biggr|\leqslant \|F_{n_j'}-F \|,\quad j\in\N,
\end{eqnarray}
we obtain  $\sup_{t\in\R}|f_{n_j'}(t)-f(t)|\to 0$ as $j\to\infty$.  Since 
\begin{eqnarray*}
	|f_{n_j'}(t)|\geqslant |f(t)|-|f_{n_j'}(t)-f(t)|\geqslant \mu -\sup_{s\in\R} |f_{n_j'}(s)-f(s)|,\quad j\in\N,
\end{eqnarray*}
it is true that $\inf_{t \in \R}|f_{n_j'}(t)|\geqslant \mu/2>0$ for all sufficiently large $j\in\N$. Hence  $F_{n_j'}\in \DS$ for these $j$, that is a contradiction.\quad $\Box$\\

\noindent
\textbf{Proof of Theorem \ref{th_Conv}.}  We will adopt here the following convention. Let $G$ be a function of bounded variation on $\R$ with characteristic function $g$, i.e. $g(t)=\int_{\R} e^{itx}\dd G(x)$, $t\in\R$. Since the uniqueness theorem is valid for functions of bounded variation, we set $\|g\|\defeq \|G\|$, where the latter is the total variation of $G$ on $\R$. So $\|g\|=0$ if and only if $g(t)=0$, $t\in\R$. It is easy to check that  $\|c\cdot g\|=|c|\cdot\|g\|$, $c\in\R$.  Let $G_1$ and $G_2$ be functions of bounded variation with characteristic functions $g_1$ and $g_2$, correspondingly. From the inequalities   $\|G_1+G_2\|\leqslant \|G_1\|+\|G_2\|$ and $\|G_1*G_2\|\leqslant \|G_1\|\cdot\|G_2\|$, we correspondingly have   $\|g_1+g_2\|\leqslant \|g_1\|+\|g_2\|$  and $\|g_1\cdot g_2\|\leqslant \|g_1\|\cdot \|g_2\|$. Moreover, such normed space of functions $g$ will be complete (see \cite{Gelfand} p. 165).

\textit{Sufficiency}. Suppose that $(i)$ and $(ii)$ hold. We define
\begin{eqnarray}\label{def_Deltan} 
	\Delta_n(t) \defeq \Ln f_n (t) -\Ln f(t),\quad t\in\R, \quad n\in\N.
\end{eqnarray}
For $n\geqslant n_0$ we have
\begin{eqnarray*}
	\Delta_n(t) &=&\sum_{u\in\la X_n\ra\setminus\{0\}}\lambda_{n,u} (e^{itu}-1)-\sum_{u\in\la X\ra\setminus\{0\}}\lambda_{u} (e^{itu}-1)\\
	&=&\sum\limits_{\substack{u\in\la X_n \ra\cup \la X\ra,\\ u\ne0}}(\lambda_{n,u}-\lambda_{u})(e^{itu}-1),\quad t\in\R.
\end{eqnarray*}
Define
\begin{eqnarray}\label{def_lambdan0lambda0}
	\lambda_{n,0}\defeq - \sum_{u\in\la X_n\ra\setminus\{0\}}\lambda_{n,u},\quad n\in\N,\quad\text{and}\quad \lambda_{0}\defeq - \sum_{u\in\la X\ra\setminus\{0\}}\lambda_{u}.
\end{eqnarray}
So we have
\begin{eqnarray}\label{eq_Deltan}
	\Delta_n(t) =\sum\limits_{u\in\la X_n \ra\cup \la X\ra}(\lambda_{n,u}-\lambda_{u})e^{itu},\quad t\in\R, \quad n\geqslant n_0.
\end{eqnarray}
Therefore 
\begin{eqnarray}\label{eq_Deltan_norm}
	\|\Delta_n\| =\sum\limits_{u\in\la X_n \ra\cup \la X\ra}|\lambda_{n,u}-\lambda_{u}|, \quad n\geqslant n_0.
\end{eqnarray}
According to \eqref{def_lambdan0lambda0}, observe that for every $n\in\N$
\begin{eqnarray*}
	|\lambda_{n,0}-\lambda_{0}|&=&\Biggl| \sum_{u\in\la X_n\ra\setminus\{0\}}\lambda_{n,u} - \sum_{u\in\la X\ra\setminus\{0\}}\lambda_{u}\Biggr| \\
	&&{}=\Biggl| \sum\limits_{\substack{u\in\la X_n \ra\cup \la X\ra,\\ u\ne0}}(\lambda_{n,u}-\lambda_{u}) \Biggr| \leqslant \sum\limits_{\substack{u\in\la X_n \ra\cup \la X\ra,\\ u\ne0}}|\lambda_{n,u}-\lambda_{u}|.
\end{eqnarray*}
Then
\begin{eqnarray*}
	\|\Delta_n\| =|\lambda_{n,0}-\lambda_{0}|+\sum\limits_{\substack{u\in\la X_n \ra\cup \la X\ra,\\ u\ne0}}|\lambda_{n,u}-\lambda_{u}| \leqslant 2\cdot\!\!\!\!\!\!\sum\limits_{\substack{u\in\la X_n \ra\cup \la X\ra,\\ u\ne0}}|\lambda_{n,u}-\lambda_{u}|, \quad n\geqslant n_0.
\end{eqnarray*}
Due to $(ii)$, we have  $\|\Delta_n\|\to 0$, $n\to\infty$. 

We next write
\begin{eqnarray*}
	f_n(t) -f(t)= f(t) \bigl(e^{\Delta_n(t)}-1\bigr)=f(t) \sum\limits_{k=1}^{\infty} \dfrac{\Delta_n(t)^k}{k!},\quad t\in\R, \quad n\geqslant n_0. 
\end{eqnarray*}
Hence
\begin{eqnarray*}
	\|F_n-F\|=\|f_n -f \|= \biggl\|f\cdot\sum\limits_{k=1}^{\infty} \dfrac{\Delta_n^k}{k!}  \biggr\|\leqslant \|f\|\cdot  \biggl\|\sum\limits_{k=1}^{\infty} \dfrac{\Delta_n^k}{k!}  \biggr\|,\quad n\geqslant n_0.
\end{eqnarray*}
 Here $\|f\|=1$ and
\begin{eqnarray*}
 \biggl\|\sum\limits_{k=1}^{\infty} \dfrac{\Delta_n^k}{k!}  \biggr\|\leqslant \sum\limits_{k=1}^{\infty} \dfrac{\|\Delta_n\|^k}{k!}= e^{\|\Delta_n\|}-1,\quad n\geqslant n_0.
\end{eqnarray*}
Thus $\|F_n-F\|\leqslant  e^{\|\Delta_n\|}-1$, $n\geqslant n_0$. Since $\|\Delta_n\|\to 0$ as $n\to\infty$, we obtain $\|F_n-F\|\to 0$,  i.e. $F_n\tovar F$.
  
\textit{Necessity}. Suppose that $F_n\tovar F$, $n\to\infty$. Then
\begin{eqnarray}\label{eq_fnfFnF_Normtozero}
	\|f_n-f\|=\|F_n-F\|\to 0,\quad n\to\infty.
\end{eqnarray}
Due to the estimate $\sup_{t\in\R}|f_n(t)-f(t)|\leqslant \|F_n-F \|$ (as in  \eqref{ineq_sup_absnorm}), $n\in\N$, 
we have that  
\begin{eqnarray}\label{conc_suptfnf}
	\sup_{t\in\R}|f_n(t)-f(t)|\to 0, \quad n\to\infty.
\end{eqnarray}
We now consider the functions $\Delta_n$, $n\in\N$, defined by \eqref{def_Deltan}:
\begin{eqnarray*}
	\Delta_n(t)=it(\gamma_n-\gamma)+\sum_{u\in\la X_n\ra\setminus\{0\}}\lambda_{n,u} (e^{itu}-1)-\sum_{u\in\la X\ra\setminus\{0\}}\lambda_{u} (e^{itu}-1),\quad t\in\R,\quad n\in\N.
\end{eqnarray*} 
Since for every $n\in\N$ the function $\Delta_n$ is continuous and $\Delta_n(0)=0$, it is the distinguished logarithm of the function $f_n/f$. Moreover, for all sufficiently large $n\in\N$
\begin{eqnarray}\label{conc_Deltan}
	\Delta_n(t)=\ln\bigl(f_n(t)/f(t)\bigr)=\ln\biggl( 1+\dfrac{f_n(t) -f(t)}{f(t)}\biggr),\quad t\in\R.
\end{eqnarray}
Indeed, by assumption,  $f\in\DS$, i.e. $f$ satisties $(S)$ with some $\mu>0$. Therefore
\begin{eqnarray*}
	\sup_{t\in\R} \biggl|\dfrac{f_n(t) -f(t)}{f(t)}\biggr| \leqslant \dfrac{\sup_{t\in\R}|f_n(t) -f(t)|}{\inf_{t\in\R}|f(t)|}\leqslant \dfrac{1}{\mu}\,\sup_{t\in\R}|f_n(t) -f(t)|,\quad n\in\N.
\end{eqnarray*}
Due to \eqref{conc_suptfnf},   the first supremum goes to zero. In particular, there exists $n_1\in\N$ such that
\begin{eqnarray}\label{ineq_fracfnf}
	\sup_{t\in\R} \biggl|\dfrac{f_n(t) -f(t)}{f(t)}\biggr| <1,\quad n\geqslant n_1.
\end{eqnarray}
Consequently, for $n\geqslant n_1$ the function in the right-hand side of \eqref{conc_Deltan} is continuous on $\R$, and it equals zero for $t=0$. By uniqueness, this function is the distinguished logarithm  of $f_n/f$, i.e.  \eqref{conc_Deltan} is valid.

We first check $(i)$. According to the above remark, from representation \eqref{def_Deltan}, on the one hand,  we see that $\sup_{t\in\R} |\Delta_n(t)|\to 0$, $n\to\infty$. On the other hand, we have the lower estimate
\begin{eqnarray*}
	\sup_{t\in\R} |\Delta_n(t)|&\geqslant&\sup_{t\in\R}|t(\gamma_n-\gamma)|-\sup_{t\in\R}\sum_{u\in\la X_n\ra\setminus\{0\}}\bigl|\lambda_{n,u} (e^{itu}-1)\bigr|-\sup_{t\in\R}\sum_{u\in\la X\ra\setminus\{0\}}\bigl|\lambda_{u} (e^{itu}-1)\bigr|\\
	&\geqslant&\sup_{t\in\R}|t(\gamma_n-\gamma)|-2\sum_{u\in\la X_n\ra\setminus\{0\}}|\lambda_{n,u}|-2\sum_{u\in\la X\ra\setminus\{0\}}|\lambda_{u}|,\quad n\in\N.
\end{eqnarray*}
Hence if there exists a strongly increasing sequence of numbers $n_k\in\N$, $k\in\N$, such that $\gamma_{n_k}-\gamma\ne 0$, $k\in\N$, then  $\sup_{t\in\R}|t(\gamma_{n_k}-\gamma)|=|\gamma_{n_k}-\gamma| \sup_{t\in\R}|t|=\infty$ and, consequently, $\sup_{t\in\R} |\Delta_{n_k}(t)|=\infty$, $k\in\N$.  In this case we have a contradiction. Therefore $\gamma_n=\gamma$ for all $n\in\N$ except a finite number, i.e. $(i)$ is valid with some $n_0\in\N$.

We now check $(ii)$. Let us fix positive integer $n\geqslant \max\{n_0,n_1\}$.
Then  \eqref{eq_Deltan} and \eqref{eq_Deltan_norm} are valid for $\Delta_n$, where we use \eqref{def_lambdan0lambda0}. Also, due to \eqref{conc_Deltan} and \eqref{ineq_fracfnf}, we have
\begin{eqnarray*}
	\Delta_n(t)=\ln\biggl( 1+\dfrac{f_n(t) -f(t)}{f(t)}\biggr)=\sum\limits_{m=1}^{\infty}\dfrac{(-1)^{m-1}}{m} \biggl( \dfrac{f_n(t) -f(t)}{f(t)}\biggr)^m,\quad t\in\R.
\end{eqnarray*}
Let us consider the function $(f_n-f)/f$. It is evident that $f_n-f$ is almost periodic function with absolutely convergent Fourier series. Since $f$ satisfies $(S)$ the function $1/f$ is also almost periodic (see \cite{Bohr} p. 35 or \cite{Levitan} p. 21)  and it has absolutely convergent Fourier series (see \cite{Gelfand} p. 175). Therefore the function $(f_n-f)/f$ is almost periodic with absolutely convergent Fourier series. For every $m\in\N$ the same is true for the function $((f_n-f)/f)^m$,  and
\begin{eqnarray*}
	\biggl\| \biggl( \dfrac{f_n -f}{f}\biggr)^m \biggr\|\leqslant \biggl\|  \dfrac{f_n -f}{f}\biggr\|^m \leqslant \bigl(\|f_n -f\| \cdot \|1/f\|\bigr)^m.
\end{eqnarray*}
Therefore
\begin{eqnarray*}
	\|\Delta_n\|= \biggl\|\sum\limits_{m=1}^{\infty}\dfrac{(-1)^{m-1}}{m} \biggl( \dfrac{f_n -f}{f}\biggr)^m\biggr\|\leqslant \sum\limits_{m=1}^{\infty}\dfrac{1}{m}\cdot \bigl(\|f_n -f\| \cdot \|1/f\|\bigr)^m.
\end{eqnarray*}
From this and \eqref{eq_fnfFnF_Normtozero} we see that $\|\Delta_n\|\to 0$ as $n\to\infty$. Thus we have
\begin{eqnarray*}
	\|\Delta_n\|=\sum\limits_{u\in\la X_n \ra\cup \la X\ra}|\lambda_{n,u}-\lambda_{u}|\to 0,\quad n\to\infty,
\end{eqnarray*}
that yields $(ii)$.\quad $\Box$\\

\noindent
\textbf{Proof of Theorem \ref{th_RelComp}.} Due to the convention \eqref{conv_lambdau_lambdanu}, we rewrite $(ii)$ in the following form
\begin{eqnarray*}
	\sup\limits_{n\in\N}\sum\limits_{k=1}^\infty|\lambda_{n,u_k}|<\infty.
\end{eqnarray*} 
This and $(iii)$ together are  equivalent to the relative compactness in $\ell^1$ of the sequence $(V_n)_{n\in\N}$, where $V_n\defeq (\lambda_{n,u_k})_{k\in\N}$, $n\in\N$, i.e. that 
every  subsequence of $(V_n)_{n\in\N}$ contains a further subsequence $(V_{n_j})_{j\in\N}$ such that $\sum_{k=1}^\infty|\lambda_{n_j,u_k} - \lambda_{u_k}| \to 0$ as $j \to \infty$, where $(\lambda_{u_k})_{k\in\N}$ is some sequence of real numbers with $\sum_{k=1}^\infty |\lambda_{u_k}|<\infty$ (see \cite{Smirnov}, p. 271 and p. 278). 

Suppose that $(i)$ holds and $(V_n)_{n\in\N}$ is relatively compact in $\ell^1$. We show that $(F_n)_{n\in\N}$ is $\DS$-relatively compact in variation.  Let us choose any subsequence $(F_{n_j})_{j\in\N}$. We first choose from $(n_j)_{j\in\N}$ a subsequence $(n_j')_{j\in\N}$ such that $(V_{n'_j})_{j\in\N}$ converges in $\ell^1$ to some sequence $V=(\lambda_{u_k})_{k\in\N}$, $k\in \N$. Next, since $\gamma_{n_j'}$, $j\in\N$, have the finite number of values from $\UnionXn$, we can choose from $(n_j')_{j\in\N}$ a subsequence $(n_j'')_{j\in\N}$ such that $\gamma_{n_j''}$ is constant for all $j\in\N$, i.e. there exists $\gamma\in \UnionXn$ such that $\gamma_{n_j''}=\gamma$, $j\in\N$.  It is easily seen that $f_{n_j''}(t)\to f(t)$, $j\to\infty$, where
\begin{eqnarray*}
	f(t)=\exp\biggl\{it\gamma+\sum_{k=1}^\infty \lambda_{u_k} (e^{itu_k}-1)\biggr\},\quad t\in\R.
\end{eqnarray*}
This function is continuous at $t=0$. By the continuity theorem, $f$ is the characteristic function for some distribution function $F$. It is discrete (because there is the absolutely convergent series in the exponent) and its characteristic function $f$ satisfies $(S)$ (see Section 2). So $F\in \DS$ and, by Theorem \ref{th_Conv}, $F_{n_j''}\tovar F$, $j\to\infty$.   
Thus we chosen from the arbitrary subsequence $(F_{n_j})_{j\in\N}$ the further subsequence $(F_{n_j''})_{j\in\N}$, which converges in variation to a distribution function from $\DS$. This yields $\DS$-relative compactness in variation for $(F_n)_{n\in\N}$.

We now suppose that $(F_n)_{n\in\N}$ is  $\DS$-relatively compact in variation. We show that $(V_n)_{n\in\N}$ is relatively compact in $\ell^1$. Let us choose any subsequence $(V_{n_j})_{j\in\N}$.  We next choose from $(n_j)_{j\in\N}$ a subsequence $(n_j')_{j\in\N}$ such that $F_{n_j'}\tovar F$, $j\to\infty$, where $F\in \DS$.
Let $f$ be its characteristic function with representation \eqref{th_Repr_eq}.  According to Theorem \ref{th_Conv}, we have
\begin{eqnarray*}
	\sum\limits_{\substack{u\in\la X_{n_j'} \ra\cup \la X\ra,\\ u\ne0}}|\lambda_{n_j',u}-\lambda_{u}|\to 0,\quad j\to\infty.
\end{eqnarray*}
Due to the convention \eqref{conv_lambdau_lambdanu}, we can write this in the following form
\begin{eqnarray*}
	\sum\limits_{\substack{u\in\UnionXn\cup \la X\ra,\\ u\ne0}}|\lambda_{n_j',u}-\lambda_{u}|\to 0,\quad j\to\infty.
\end{eqnarray*}
It is clear that if $\la X\ra \setminus \UnionXn\ne \varnothing$ then $\lambda_{u, n_j'}=0$, and hence $\lambda_u=0$ for all $u\in\la X\ra \setminus \UnionXn$. Therefore we have   $\sum_{u\in\UnionXn\setminus\{0\}}|\lambda_{n_j',u}-\lambda_{u}|=\sum_{k=1}^\infty|\lambda_{n_j',u_k}-\lambda_{u_k}|\to 0$, $j\to\infty$. This proves the relative compactness $(V_n)_{n\in\N}$ in $\ell^1$.

We now show that $(i)$ holds. Suppose, contrary to our claim, that there exists subsequence $(\gamma_{n_j})_{j\in\N}$ with distinct values.  We choose from $(n_j)_{j\in\N}$ a subsequence $(n_j'')_{j\in\N}$ such that $F_{n_j''}\tovar F$, $j\to\infty$ with some $F\in \DS$.  According to Theorem \ref{th_Conv},  $\gamma_{n_j''}$ is constant for all sufficiently large $j\in\N$, a contradiction. \quad $\Box$\\

\noindent
\textbf{Proof of Theorem \ref{th_StochComp}.}  \textit{Sufficiency}. Let us choose any  subsequence $(F_{n_j})_{j\in\N}$.  Since $(F_n)_{n\in\N}$ is $\DS$-relatively compact in variation, we can choose futher  subsequence $(F_{n_j'})_{j\in\N}$, which converges in variation to some distribution function $F\in \DS$. Let $f$ be its characteristic function with representation \eqref{th_Repr_eq}.  If $F$ is degenerate, then  $f(t)=e^{it\gamma}$, i.e. $\lambda_{u}=0$, $u\in \la X\ra\setminus\{0\}$. By Theorem \ref{th_Conv}, we have that $\gamma_{n_j'}=\gamma$ for all sufficiently large $j\in\N$, and 
\begin{eqnarray*}
	\sum_{\substack{u\in\la X_{n_j'} \ra \cup \la X\ra\\ u\ne 0}}|\lambda_{n_j',u}|\to 0,\quad j\to\infty.
\end{eqnarray*}
This contradicts \eqref{th_StochComp_cond}. Hence $F$ is always nondegenerate distribution function. Thus we have $\DS$-stochastic compactness in variation of $(F_n)_{n\in\N}$. 

\textit{Necessity}. Suppose that $(F_n)_{n\in\N}$ is $\DS$-stochastically compact in variation. Then, by definition, it is $\DS$-relatively compact in variation. Suppose that \eqref{th_StochComp_cond} is not true, i.e. there exists a sequence $(n_j)_{j\in\N}$ such that $\sum_{u\in\la X_{n_j} \ra\setminus\{0\}}|\lambda_{n_j,u}|\to 0$, $j\to\infty$. By Theorem \ref{th_RelComp}, $\gamma_{n_j}$, $j\in\N$, have finite number of values. Let us choose from $(n_j)_{j\in\N}$ a subsequence $(n_j')_{j\in\N}$    such that $\gamma_{n_j'}=\gamma$, $j\in\N$, with some $\gamma\in\R$. Then, according to Theorem \ref{th_Conv}, $F_{n_j'}\tovar F$, $j\to\infty$, where $F$ is degenerate distribution function with unique growth point $\gamma$. This contradicts to $\DS$-stochastic compactness in variation for $(F_n)_{n\in\N}$. \quad $\Box$

\section{Acknowledgments}
The work of I. A. Alexeev was supported by Ministry of Science and Higher Education of the Russian Federation, agreement  075-15-2019-1619.  The work of A. A. Khartov was supported by RFBR--DFG grant 20-51-12004.

\end{document}